\documentclass{article}
\usepackage{amsmath}
\usepackage{amsfonts}

    \newtheorem{theorem}{Theorem}[section]
    \newtheorem{lemma}[theorem]{Lemma}
    \newtheorem{proposition}[theorem]{Proposition}

\newenvironment{proof}[1][Proof]{\begin{trivlist}
\item[\hskip \labelsep {\bfseries #1}]}{\end{trivlist}}
\newcommand{\qed}{\nobreak \ifvmode \relax \else
      \ifdim\lastskip<1.5em \hskip-\lastskip
      \hskip1.5em plus0em minus0.5em \fi \nobreak
      \vrule height0.75em width0.5em depth0.25em\fi}

\begin{document}
\begin{title}{The Expansion for the overlap function}
\end{title}
\author{S\'ergio de Carvalho Bezerra\\
\it{Instituto de Matem\'atica e Estat\'istica-USP-Brazil}\\sergio@ime.usp.br}
\maketitle
\begin{abstract} In this work, it is proved the complete expansion for the
  second moment of the overlap function for the Sherrington-Kirkpatrick
  model. It is a technical result which takes advantage of the cavity method
  and other induction arguments.
\end{abstract}
\footnotetext{The author has been supported by the FAPESP.}
\section{Introduction}
\subsection{SK model}
This work concernes the study of a quantity that plays an important role in
the Sherrington-Kirkpatrick model (SK model). In the following, I describe the
SK model. In this model, we have the particularity to obtain two levels of
randomness. More precisely, we suppose that a certain system has $N$-particles and the set of all possible system configurations
 is $\Sigma_N=\{-1,+1\}^N$. We define the probability of realization of one
 configuration $\sigma\in\Sigma_N$ as:
$$
\mathrm P_N(\sigma)=\frac{e^{\beta \mathrm E(\sigma)}}{Z_N},
$$   
where $\beta$ is a parameter of the system (a positive real number), the term
$Z_N=\sum_{\sigma\in\Sigma^N}e^{\beta \mathrm E(\sigma)}$ is a normalization factor
and the term $\mathrm E(\sigma)$ represents the energy of the configuration
$\sigma$. In this way, we obtain our first probabilistic model
$(\Sigma_N,\mathcal{F}_N,\mathrm P_N)$. However, the measure $\mathrm P_N$ is a random measure
because the energy of the configuration $\sigma=(\sigma_1,\ldots,\sigma_N)$ is
a random variable defined by:
$$\mathrm E(\sigma)=\frac{1}{\sqrt{N}}\sum_{1\leq i<j\leq
  N}\sigma_i\sigma_jg_{i,j}.$$
The $(g_{i,j})_{1\leq i<j\leq N}$ is a family of independent standard gaussian
  random variables that live in a second probabilistic world
  $(\Omega_N,\mathcal{G}_N,\mathrm G_N)$. When we look at the correlation
  between the energy levels of two different configurations $\sigma^1$ and
  $\sigma^2$, we obtain that
$$
\mathbf{E}\left[\mathrm E(\sigma^1)\mathrm E(\sigma^2)\right]=\frac{1}{N}\sum_{1\leq
  i<j\leq N}\sigma^1_i\sigma^1_j\sigma^2_i\sigma^2_j= \frac{1}{2N}\left[\left(\sum_{ i=1}^N\sigma^1_i\sigma^2_i\right)^2-N\right].
$$
Consequently, the correlation depends on the quantity
$$
R_{1,2}=\frac{1}{N}\sum_{i=1}^N\sigma_i^1\sigma^2_i
$$
 called overlap function. 

If we have a function $f$ defined in the product space $\Sigma_N^n$, we will
denothe by $\langle f\rangle$ the expectation with respect to the product
measure $P_N^{\otimes n}$ and by $\nu(f)$ the double expectation $\mathbf E
\langle f\rangle$.

\subsection{Previous Results}

 According to the literature, we can obtain result for two different regions:
 the one of high temperature (small $\beta$) and the one of low temperature
 (large $\beta$). The book \cite{Ta} is a good survey of the results for the
 SK model and its generalizations from a mathematical point of view. During this work, we will be focus on the high temperature region.
 We can find, in the book \cite{Ta}, the following result:
\begin{theorem}Given $\beta<1$, we have that:
\begin{align}
\nu(R^2_{1,2})=\frac{1}{N(1-\beta^2)}+\mathrm O\left(\frac{1}{N^{3\slash2}}\right).\nonumber
\end{align}
\end{theorem}

In the article~\cite{BMRT}, we can see the result.
\begin{theorem}Given $\beta<1$, we have that:
\begin{align}
\nu(R_{1,2}^2)=\frac{1}{N(1-\beta^2)}-\frac{\beta^2(1+\beta^2)}{N^2(1-\beta^2)^4}+\mathrm O\left(\frac{1}{N^{\frac{5}{2}}}\right),\nonumber
\end{align}
\end{theorem} 
and they also obtained numerical results for the third term of the expansion.

Our goal is to obtain the complete expansion for $\nu(R_{1,2}^2)$.
\subsection{Cavity Method}

At this point, we introduce a family of fonctionals
$\nu_t$ that depends on a parameter $t\in [0,1]$. First of all, given one function
$f:\Sigma_N^n\rightarrow \mathbb{R}$, we define $\nu_t(f)=\mathbf{E}\langle
f\rangle_t$ where:
\begin{align}
\langle f\rangle_t=\sum_{\sigma^1,\ldots,\sigma^n\in
  \Sigma_N}\frac{f(\sigma^1,\ldots,\sigma^n)e^{\beta(\mathrm
  E_t(\sigma^1)+\ldots+\mathrm E_t(\sigma^n))}}{Z_{N,t}^n}\nonumber
\end{align}
with the normalization factor given by:
\begin{align}
 Z_{N,t}=\displaystyle\sum_{\sigma \in\Sigma_N}e^{\beta\mathrm
 E_t(\sigma)}
\nonumber
\end{align}
and the energy by:
$$
\mathrm E_t(\sigma)=\frac{1}{\sqrt{N}}\sum_{1\leq i<j\leq
  N-1}\sigma_i\sigma_jg_{i,j}+\left(\frac{t}{N}\right)^{\frac{1}{2}}\sigma_N\sum_{1\leq
i\leq N-1}\sigma_ig_{i,N}.
$$
Therefore, if the parameter $t$ takes the value one we have the usual energy
of one configuration with $N$ particles else $t$ assumes the value zero we have the energy of one
configuration with $N-1$ particles for a different temperature. Thus, we see the idea of induction with
the introduction of the parameter $t$.

\subsection{Some notations}

The first idea used in this work is to perform a Taylor`s series for the function $t\rightarrow
\nu_t(f)$. When we do it the terms $\epsilon_{l_1}=\sigma_N^{l_1}$ appear for $l_1$ a
positive integer and also appear the functions
$R^-_{l_1,l_1^{\prime}}=\frac{1}{N}\sum_{1\leq i\leq
  N-1}\sigma_i^{l_1}\sigma_i^{l_1^{\prime}}$ for $l_1^{\prime}$ another positive
integer greater than $l_1$. The first differentiation of $\nu_t(f)$ satisfies:
\begin{proposition}\label{P:derivada} Given a function $f$ on $\Sigma_N^n$ and $t\geq 0$, we
  have
\begin{align}
\nu_t^{\prime}(f)&=\beta^2\sum_{1\leq l_1<l_1^{\prime}\leq
  n}\nu_t(f\epsilon_{l_1}\epsilon_{l_1^{\prime}}R^-_{l_1,l_1^{\prime}})\nonumber\\
&-\beta^2n\sum_{l_1\leq n}\nu_t(f\epsilon_{l_1}\epsilon_{n+1}R_{l_1,n+1}^-)\nonumber\\
&+\beta^2\frac{n(n+1)}{2}\nu_t(f\epsilon_{n+1}\epsilon_{n+2}R_{n+1,n+2}^-)\nonumber
\end{align}
\end{proposition}
which the proof can be found in \cite{Ta}. Let $\ell$ be an integer greater than
or equal to one. When we take the $\ell^{th}$-differentiation of $\nu_t(f)$ it appears the terms: 
\begin{align}
S^-_{l}=\prod_{i=1}^{\ell}\epsilon_{l_i}\epsilon_{l_i^{\prime}}R^-_{l_i,l_i^{\prime}},\nonumber
\end{align}
with $l=(l_1,l_1^{\prime},\ldots,l_{\ell},l_{\ell}^{\prime})$.
These terms define naturally $3^{\ell}$ sets of $sequences$ of the type
$(l_1,l^{\prime}_1,\ldots,l_{\ell},l_{\ell}^{\prime})$ denoted by $C^f_{\ell,m}$ with
$m=1,\ldots,3^{\ell}$ and defined for $\ell=1$ by: 
\begin{align}
C^f_{1,1}&=\{(l_1,l_1^{\prime})|1\leq l_1<l_1^{\prime}\leq
n\}\nonumber\\
C^f_{1,2}&=\{(l_1,l_1^{\prime})|1\leq l_1\leq n;\;l_1^{\prime}=n+1\}\nonumber\\
C^f_{1,3}&=\{(l_1,l_1^{\prime})|l_1=n+1;\;l_1^{\prime}=n+2\}\nonumber
\end{align}
where $n$ is the number of configurations that $f$ depends on. For $\ell\geq 2$ and each positive integer of the set $\left\{\lceil
 \frac{m}{3}\rceil|m\in\{1,\ldots,3^{\ell}\}\right\}$ we define recursively
 the numbers
\begin{align}\label{D:maximo}
M^f_{\ell-1,\lceil \frac{m}{3}\rceil}=\{\mbox{ the number of configurations of
}fS_{\hat l}^-\mbox{ with }\hat l\in C_{\ell-1,\lceil \frac{m}{3}\rceil}\}
\end{align} 
and the sets:
\begin{align}
C^f_{\ell,3\lceil \frac{m}{3}\rceil-2}&=\{(l_1,\ldots,l_{\ell},l_{\ell}^{\prime})|(l_1,\ldots,l_{\ell-1}^{\prime})\in
C^f_{\ell-1,\lceil \frac{m}{3}\rceil},1\leq l_{\ell}<l_{\ell}^{\prime}\leq
M^f_{\ell-1,\lceil \frac{m}{3}\rceil}\},\nonumber\\
C^f_{\ell,3\lceil\frac{m}{3}\rceil-1}&=\{(l_1,\ldots,l_{\ell},l_{\ell}^{\prime})|(l_1,\ldots,l_{\ell-1}^{\prime})\in
C^f_{\ell-1,\lceil \frac{m}{3}\rceil},1\leq l_{\ell}\leq M^f_{\ell-1,\lceil \frac{m}{3}\rceil}\mbox{ and
}\nonumber\\&l_{\ell}^{\prime}=M^f_{\ell-1,\lceil \frac{m}{3}\rceil}+1\},\nonumber\\
C^f_{\ell,3\lceil \frac{m}{3}\rceil}&=\{(l_1,\ldots,l_{\ell},l_{\ell}^{\prime})|(l_1,\ldots,l_{\ell-1}^{\prime})\in
C^f_{\ell-1,\lceil \frac{m}{3}\rceil},l_{\ell}=M^f_{\ell-1,\lceil \frac{m}{3}\rceil}+1\mbox{ and }\nonumber\\
&l_{\ell}^{\prime}=M^f_{\ell-1,\lceil \frac{m}{3}\rceil}+2\}.
\end{align}

We observe that the numbers $M^f_{\ell,m}$ and the sets
$C^f_{\ell,m}$ are well defined. For $\ell=1$ we have that
$M_{1,1}=n,M_{1,2}=n+1$ and $M_{1,3}=n+2$. By hypothesis of induction on $\ell$
we suppose that the numbers $M^f_{\ell-1,\lceil\frac{m}{3}\rceil}$ and the
sets $C^f_{\ell-1,\lceil\frac{m}{3}\rceil}$ are well defined. Consequently,
for each $m\in\{1,\ldots,3^{\ell}\}$ by their definitions the sets
$C^f_{\ell,3\lceil\frac{m}{3}\rceil-2},C^f_{\ell,3\lceil\frac{m}{3}\rceil-1}$
and $C^f_{3\lceil\frac{m}{3}\rceil}$ are also well defined and the numbers
$M^f_{\ell,3\lceil\frac{m}{3}\rceil-2}=M^f_{\ell-1,\lceil\frac{m}{3}\rceil},M^f_{\ell,3\lceil\frac{m}{3}\rceil-1}=M^f_{\ell-1,\lceil\frac{m}{3}\rceil}+1$
and
$M^f_{\ell,3\lceil\frac{m}{3}\rceil}=M^f_{\ell-1,\lceil\frac{m}{3}\rceil}+2$.

Now, we also introduce the following notation: let $u$ be a integer between $1$
and $\ell$. Given any sequence
$l=(l_1,l_1^{\prime},\ldots,l_{\ell},l_{\ell}^{\prime})$ we define for each $u$ the set 
\begin{align}
  G_l^u=\left\{l^u=(l_{i_1},l_{i_1}^{\prime},\ldots,l_{i_u},l_{i_u}^{\prime})|
  (i_1,\ldots,i_u)\in\{1,\ldots,\ell\}\mbox{ and are all different}\right\}.\nonumber
\end{align}
and given a subsequence $l^u=(l_{i_1},l_{i_1}^{\prime},\ldots,l_{i_u},l_{i_u}^{\prime})\in G_l^u$ we also define the term
$\epsilon^{l^u}=\prod_{r=1}^u\epsilon_{l_{i_r}}\epsilon_{l_{i_r}^{\prime}}$,
the term
$\epsilon^{\hat l^u}=\epsilon^l\epsilon^{l^u}$ and the functions
\begin{align}
R_{l^u}=R_{l_{i_1},l_{i_1}^{\prime}}\ldots R_{l_{i_u},l_{i_u}^{\prime}} \mbox{
and }R_{l^u}^-=R_{l_{i_1},l_{i_1}^{\prime}}^-\ldots R_{l_{i_1},l_{i_1}^{\prime}}^-.\nonumber
\end{align}

We also introduce the coefficient $\rho^f(l)$ defined in the Proposition below
which the proof can be found in the section \textbf{A Result of Differentiation}:
\begin{proposition}\label{P:derivada-l}Given a function $f$ on $\Sigma_N^n$ and $t\geq0$, we
  obtain that for a positive integer $\ell$:
\begin{align}
\nu^{(\ell)}_t(f)=\sum_{m=1}^{3^{\ell}}\sum_{l\in C^f_{\ell,m}}\beta^{2\ell}\rho^f(l)\nu_t(fS^-_l),\nonumber
\end{align}
 The
coefficients $\rho^f(l)$ are
defined for $\ell$ equals to one by: $\rho^f(l)=1$ if $l\in C^f_{1,1}$,
$\rho^f(l)=-n$ if $l\in C^f_{1,2}$ and $\rho^f(l)=\frac{n(n+1)}{2}$ if $l\in
C^f_{1,3}$. Given $\ell$ greater or equal to two they are defined recursively
by: we consider that $l=(\hat l,l_{\ell},l_{\ell}^{\prime})$ with $\hat
l\in C^f_{\ell-1,\lceil \frac{m}{3}\rceil}$ for some $m$
then
\begin{align}
\rho^f(l)=\begin{cases}\rho^f(\hat l),&\mbox{ if } l\in
C^f_{\ell,3\lceil\frac{m}{3}\rceil-2},\\
-\rho^f(\hat l)M^f_{\ell-1,\lceil\frac{m}{3}\rceil},& \mbox{ if }l\in
C^f_{\ell,3\lceil\frac{m}{3}\rceil-1}\mbox{ and }\\
\rho^f(\hat
l)\frac{M^f_{\ell-1,\lceil\frac{m}{3}\rceil}(M^f_{\ell-1,\lceil\frac{m}{3}\rceil}+1)}{2} &
\mbox{ if }l\in C^f_{\ell,3\lceil\frac{m}{3}\rceil}
\end{cases}\nonumber
\end{align}
\end{proposition}

\subsection{Contributions}

The main results of this work are the following ones:
\begin{lemma}\label{L:main}Given $\beta$ a positive real number less than one. We determine for all
  sequence  $l=(l_1,l_1,\ldots,l_{\ell},l_{\ell}^{\prime})$ such that
  $l_i<l_i^{\prime}$ for all $i=1,\ldots,\ell$ and
  $\epsilon^l=1$ with $\ell\leq 2k$ for a positive integer $k$ and $l=(\hat
  l,l_{\ell},l_{\ell}^{\prime})$, the coefficients
  $\lambda_j^l$ such that
\begin{align}
\nu_0(\epsilon^lR_l^-)=\sum_{j=0}^k\frac{\lambda_j^l}{N^{j}}+\mathrm O\left(\frac{1}{N^{k+\frac{1}{2}}}\right).\nonumber
\end{align}
\begin{itemize}
\item[(a)] if $k$ is equal to one the only possible value for $\ell$ is $2$ and we
  have
\begin{align}
\lambda_1^l=\frac{1}{(1-\beta^2)},\nonumber
\end{align}  
\item[(b)]if $k\geq 2$ and $\ell=2k$ it enoughs to know the value of the coefficients
  $\lambda_j^{\gamma}$ for $j=1,\ldots,k-1$ of all sequences
  $\gamma=(\gamma_1,\gamma_1^{\prime},\ldots,\gamma_{z},\gamma_{z}^{\prime})$ such that
  $z\leq 2k-2$ and
  $\epsilon^{\gamma}=1$. Then,
\begin{align}
\lambda_k^l&=\frac{1}{1-\beta^2}\Bigg[\sum_{u=1}^{2k-2}\sum_{l^u\in G_l^u}\sum_{v=1}^u\sum_{l^{u,v}\in
  G_{l^u}^v}(-1)^{2k-u}\lambda_{v-k}^{l^{u,v}}+\sum_{u=1}^{2k-2}\sum_{\hat
l^u\in G_{\hat l}^u}\lambda_{u+1-k}^{\hat
  l^u}+\nonumber\\&+\sum_{u=1}^{2k-2}\sum_{\hat l^u\in G_{\hat
  l}^u}\sum_{r=1}^{2k-2-u}\sum_{p=1}^{3^r}\sum_{\eta\in
  C_{r,p}^{\epsilon^{\hat l^u}R_{\hat
  l^u}}}\beta^{2r}\rho^{\epsilon^{\hat l^u}R_{\hat l^u}}(\eta)\lambda_{u+1-k}^{(\hat
  l^u,\eta)}\Bigg]\nonumber
\end{align}
and the coefficients $\lambda_j$ with $j$ negative are all of them equal to zero.
\item[(c)]if $k\geq 2$ and $2\leq \ell=\hat k< 2k$ it enoughs to know the
  coefficients $\lambda_j^{\gamma}$ for $j=1,\ldots,k-1$ for all sequences
  $\gamma=(\gamma_1,\gamma_1^{\prime},\ldots,\gamma_z,\gamma_z^{\prime})$ such
  that $z\leq 2k-2$ and $\epsilon^{\gamma}=1$ and $\lambda_k^{\gamma}$ if $\hat
  k<z\leq 2k$ and $\epsilon^{\gamma}=1$. Thus, we obtain
\begin{align}
\lambda_k^l&=\frac{1}{1-\beta^2}\Bigg[\sum_{u=1}^{\hat k-1}\sum_{l^u\in
  G_l^u}\sum_{v=1}^{u}\sum_{l^{u,v}\in G_{l^u}^v}(-1)^{\hat k-u}\lambda_{k-\hat
  k+v}^{l^{u,v}}+\mathrm{I}_{\{\hat k=1+k\}}\nonumber\\&+\sum_{u=1}^{\hat k-2}\sum_{\hat l^u\in G_{\hat
  l}^u}\Big[\lambda_{k-\hat k+1+u}^{\hat
  l^u}+\sum_{r=1}^{2k-2-u}\sum_{p=1}^{3^r}\sum_{\begin{array}{c}\eta\in
  C_{r,p}^{\epsilon^{\hat l}R_{\hat l}}\\\epsilon^{\hat
  l^u}\epsilon^{\eta}=1\end{array}}\frac{\beta^{2r}\rho^{\epsilon^{\hat
  l}R_{\hat l}}(\eta)}{r!}\lambda_{k-\hat
  k+1+u}^{(\hat
  l^u,\eta)}\Big]\nonumber\\&\sum_{r=2}^{2k-\hat k+1}\sum_{p=1}^{3^r}\sum_{\begin{array}{c}\eta\in
  C_{r,p}^{\epsilon^{\hat l}R_{\hat l}}\\\epsilon^{\hat
  l}\epsilon^{\eta}=1\end{array}}\frac{\beta^{2r}\rho^{\epsilon^{\hat
  l}R_{\hat l}}(\eta)}{r!}\lambda_k^{(\hat
  l,\eta)}
-\sum_{r=1}^{2k-\hat k}\sum_{p=1}^{3^r}\sum_{\begin{array}{c}\eta\in
  C_{r,p}^{\epsilon^lR_l}\\\epsilon^{\eta}=1\end{array}}\Big(\lambda_{k-\hat k}^{\eta}\nonumber\\&\frac{\beta^{2r}\rho^{\epsilon^lR_l}(\eta)}{r!}\Big)-\sum_{r=1}^{2k-\hat k}\sum_{p=1}^{3^r}\sum_{\eta\in
  C_{r,p}^{\epsilon^lR_l}}\frac{\beta^{2r}\rho^{\epsilon^lR_l}(\eta)}{r!}\sum_{u=1}^{\hat
  k-1}\sum_{\begin{array}{c}l^u\in G_l^u\\\epsilon^{l^u}\epsilon^{\eta}=1
  \end{array}}\lambda_{k-\hat k+u}^{(l^u,\eta)}\Bigg],\nonumber
\end{align}
where $\mathrm I_{\hat k=1+k}$ is the indicator function of the set $\{\hat
k|\hat k=k+1\}$.
\end{itemize}
\end{lemma}

In this work, we also demonstrate the complete expansion for the
Expression of $\nu(R_{1,2}^2)$ :
\begin{theorem}\label{T:main} Let $\beta$ be less than one and $k_0$ be an
  integer positive $k_0\geq 1$. Then,
\begin{align}
\nu(R_{1,2}^2)=\frac{1}{N}+\sum_{j=1}^{k_0}\frac{c_j(\beta)}{N^j}+\mathrm
O\left(\frac{1}{N^{k_0+\frac{1}{2}}}\right)\nonumber
\end{align}
with
$$c_{j}(\beta)=\sum_{\ell=1}^{2k_0-1}\sum_{m=1}^{3^r}\sum_{\begin{array}{c}l\in
  C_{\ell,m}^{\epsilon_1\epsilon_2R_{1,2}^-}\\\epsilon_1\epsilon_2\epsilon^l=1\end{array}}\frac{\beta^{2\ell}\rho^{\epsilon_1\epsilon_2R_{1,2}^-}(l)}{\ell!}\lambda_j^{(1,2,l)}$$
and the coefficients $\lambda_j^{(1,2,l)}$ are defined in the
Lemma~\ref{L:main}.
\end{theorem}

This Theorem is quoted by Talagrand in \cite{Ta} as a research problem.

In the next section, we prove the Proposition~\ref{P:derivada-l}. In the third
section, we obtain some preliminary results. In the following section, we
prove the Lemma~\ref{L:main} and in the last section, we obtain the proof of Theorem~\ref{T:main}.

\section{A Result of Differentiation}

Our first result concernes the differentiation of order $\ell$ for a function $f$ which produces:

\begin{proof}[The proof of Proposition~\ref{P:derivada-l}]The proof is by induction on $\ell$. The case $\ell$ equals to
  one is the result obtained in the
proposition~\ref{P:derivada}. We suppose that the result is true for
$\ell-1$. Then, we have that
\begin{align}
\nu_t^{\ell-1}(f)=\sum_{m=1}^{3^{\ell-1}}\sum_{\hat l\in
  C^f_{\ell-1,m}}\beta^{2(\ell-1)}\rho^f(\hat l)\nu_t(fS_{\hat l}^-)\nonumber
\end{align}
At this point, we remember the definition of the numbers $M^f_{\ell,m}$ given by
Equation~(\ref{D:maximo}) and we
apply the Propositon~\ref{P:derivada} for the function
$g=fS_{\hat l}$, then we obtain that

\begin{align}
\nu_t^{\ell}(f)=\sum_{m=1}^{3^{\ell-1}}&\sum_{\hat l\in
  C^f_{\ell-1,m}}\beta^{2(\ell-1)}\rho^f(\hat l)\Bigg(\sum_{1\leq
  l_{\ell}<l_{\ell}^{\prime}\leq
  M^f_{\ell-1,m}}\beta^2\nu_t(fS_l^-S_{l_{\ell},l_{\ell}^{\prime}}^-)+
  \nonumber\\
&-\sum_{1\leq l_{\ell}\leq M^f_{\ell-1,m}}\beta^2M^f_{\ell-1,m}\nu_t(fS_{\hat
  l}^-S_{l_{\ell},M^f_{\ell-1,m}}^-)\nonumber\\
&+\beta^2\frac{M^f_{\ell-1,m}(M^f_{\ell-1,m}+1)}{2}\nu_t(fS_{\hat l}^-S_{M^f_{\ell-1,m}+1,M^f_{\ell-1,m}+2}^-)\Bigg)\nonumber
\end{align}
 and the result follows when we look at the definition of the sets
 $C^f_{\ell,m}$ and the fact that $S^-_{\hat
 l}S^-_{l_{\ell},l_{\ell}^{\prime}}=S_{\hat l,l_{\ell},l_{\ell}^{\prime}}^-$\qed
\end{proof} 

Using the Proposition~\ref{P:derivada-l} and perfoming a Taylor series, we
obtain that
\begin{align}\label{E:taylor-f}
\nu(f)=\nu_0(f)+\sum_{\ell=1}^{\infty}\sum_{m=1}^{3^{\ell}}\sum_{l\in C^f_{\ell,m}}\frac{\beta^{2\ell}\rho^f(l)}{\ell!}\nu_0(fS^-_l)
\end{align}

Some terms in the above expression vanish, in reason of the following result:

\begin{proposition}\label{P:media-nula} Let $f$ be a function defined on $\Sigma_N^n$. Assume
  $f=f^-f^{\prime}$ where $f^-$ is a function of the $N-1$-system, and
  $f^{\prime}$ depends only on $\epsilon_1,\ldots,\epsilon_n$. If $Av\,
  f^{\prime}=0$ (where $Av$ means average on $\epsilon_1=\pm
  ,\ldots,\epsilon_n=\pm 1$) then
\begin{align}
\nu_0(f)=0.\nonumber
\end{align} 
\end{proposition}
 The proof of this Proposition can be also found in \cite{Ta}.

\section{The case $\epsilon^lR_l$}

We can write for any sequence $l=(l_1,l_1^{\prime},\ldots,l_{\ell},l_{\ell}^{\prime})$ that 
\begin{align}\label{E:R-to-R^-}
\epsilon^lR_l&=\epsilon^l\prod_{i=1}^{\ell}\left(R_{l_i,l_i^{\prime}}^-+\frac{\epsilon_{l_i}\epsilon_{l_i^{\prime}}}{N}\right)=\epsilon^l\left(\frac{\epsilon^l}{N^{\ell}}+\sum_{u=1}^{\ell}\sum_{l^u\in
      G_l^u}\frac{R_{l^u}^-\epsilon^{\hat
      l^u}}{N^{\ell-u}}\right)\nonumber\\
&=\frac{1}{N^{\ell}}+\sum_{u=1}^{\ell}\sum_{l^u\in G_l^u}\frac{\epsilon^{l^u}R^-_{l^u}}{N^{\ell-u}}
\end{align}
which produces the following result:
\begin{proposition}\label{P:nu-R} Given any sequence
  $l=(l_1,l_1^{\prime},\ldots,l_{\ell},l_{\ell}^{\prime})$ such that
  $l_i<l_i^{\prime}$ for all $i=1,\ldots,\ell$, we obtain that
\begin{align}
\nu(\epsilon^lR_l)&=\frac{1}{N^{\ell}}+\sum_{u=1}^{\ell}\sum_{l^u\in
    G_l^u}\frac{\nu_0(\epsilon^{l^u}R_{l^u}^-)}{N^{\ell-u}}\nonumber\\&+\sum_{u=1}^{\ell}\sum_{l^u\in
    G_l^u}\sum_{r=1}^{\infty}\sum_{p=1}^{3^r}\sum_{\begin{array}{c}\eta\in
    C^{\epsilon^{l^u}R_{l^u}^-}_{r,p}\\
\epsilon^{l^u}\epsilon^{\eta}=1\end{array}}\frac{1}{N^{\ell-u}}\frac{\beta^{2r}\rho^{\epsilon^{l^u}R_{l^u}^-}(\eta)}{r!}\nu_0(R_{l^u}^-R_{\eta}^-)\nonumber
\end{align}
\end{proposition}
\begin{proof}
We start by applying the Equation~(\ref{E:taylor-f}) for each one of the
functions in Equation~(\ref{E:R-to-R^-}), in reason of the linearity of $\nu$, then we obtain that
\begin{align}
\nu(\epsilon^lR_l)&=\frac{1}{N^{\ell}}+\sum_{u=1}^{\ell}\sum_{l^u\in
G_l^u}\frac{1}{N^{\ell-u}}\Bigg(\nu_0(\epsilon^{l^u}R_{l^u}^-)+\nonumber\\&+\sum_{r=1}^{\infty}\sum_{p=1}^{3^r}\sum_{\eta\in
C^{\epsilon^{l^u}R_{l_u}^-}_{r,p}}\frac{\beta^{2r}\rho^{\epsilon^{l^u}R_{l^u}^-}(\eta)}{r!}\nu_0(\epsilon^{l^u}R_{l^u}^-S_{\eta}^-)\Bigg)\nonumber
\end{align}

We remark that $S_{\eta}^-=\epsilon^{\eta}R_{\eta}^-$. Consequently,  We can
make the constraint that $\epsilon^{\eta}\epsilon^{l^u}=1$ in reason of the
Proposition~\ref{P:media-nula} which concludes the proof\qed 
\end{proof}

We also prove another relation:
\begin{proposition}\label{P:nu_0-R-} Let $l$ be a sequence
  $(l_1,l_1^{\prime},\ldots,l_{\ell},l_{\ell}^{\prime})$ and $\hat
  l=(l_1,\ldots,l^{\prime}_{\ell-1})$ such that $\epsilon^l=1$ and
  $l_i<l_i^{\prime}$ for all $i=1,\ldots,\ell$. Then, we obtain that
\begin{align}
&\nu_0(\epsilon^lR_l^-)=\frac{(-1)^{\ell}}{N^{\ell}}+\sum_{u=1}^{\ell-1}\sum_{l^u\in
G_l^u}\frac{(-1)^{\ell-u}}{N^{\ell-u}}\nu(\epsilon^{l^u}R_{l^u})+\nu(\epsilon^{\hat
l}R_{\hat l})\nonumber\\
&\sum_{u=1}^{\ell}\sum_{l^u\in
  G_l^u}\sum_{r=1}^{\infty}\sum_{p=1}^{3^r}\sum_{\begin{array}{c}\eta\in
    C^{\epsilon^{l^u}R_{l^u}}_{r,p}\\\epsilon^{\eta}=1\end{array}}\frac{(-1)^{\ell-u+1}\beta^{2r}\rho^{\epsilon^{l^u}R_{l^u}}(\eta)}{r!N^{\ell}}\nu_0(\epsilon^{\eta}R_{\eta}^-)+
\sum_{u=1}^{\ell}\nonumber\\&\sum_{l^u\in
G_l^u}\sum^{\infty}_{r=1}\sum_{p=1}^{3^r}\sum_{\eta\in
  C^{\epsilon^{l^u}R_{l^u}}_{r,p}}\sum_{v=1}^{u}\sum_{\begin{array}{c}l^{u,v}\in
G_{l^u}^v\\\epsilon^{l^{u,v}}\epsilon^{\eta}=1\end{array}}\frac{(-1)^{\ell-u+1}\beta^{2r}\rho^{\epsilon^{l^u}R_{l^u}}(\eta)}{r!N^{\ell-v}}\nu_0(R_{l^{u,v}}^-R_{\eta}^-)\nonumber
\end{align}
we observe that if $\epsilon^l\neq 1$ the Equation is also true but it is equal to zero.
\end{proposition}
\begin{proof}
Given the sequence $l=(l_1,l_1^{\prime},\ldots,l_{\ell},l_{\ell}^{\prime})$ we
have that:
\begin{align}
\epsilon^lR_l^-&=\epsilon^l\prod_{i=1}^{\ell}\left(R_{l_i,l_i^{\prime}}-\frac{\epsilon_{l_i}\epsilon_{l_i^{\prime}}}{N}\right)\nonumber\\
&=\frac{(-1)^{\ell}}{N^{\ell}}+\sum_{u=1}^{\ell}\sum_{l^u\in G^u_l}\frac{(-1)^{\ell-u}\epsilon^{l^u}R_{l^u}}{N^{\ell-u}}\nonumber
\end{align}
Thus, we obtain by linearity that
\begin{align}
\nu_0(\epsilon^lR_l^-)=\frac{(-1)^{\ell}}{N^{\ell}}+\sum_{u=1}^{\ell}\sum_{l^u\in
G_l^u}\frac{(-1)^{\ell-u}}{N^{\ell-u}}\nu_0(\epsilon^{l^u}R_{l^u})\nonumber
\end{align}
When, we apply above the Equation~(\ref{E:taylor-f}) for each one of the term in
the right side, we get
\begin{align}
\nu_0(\epsilon^lR_l^-)&=\frac{(-1)^{\ell}}{N^{\ell}}+\sum_{u=1}^{\ell}\sum_{l^u\in
G_{l}^u}\frac{(-1)^{\ell-u}}{N^{\ell-u}}\nu(\epsilon^{l^u}R_{l^u})\nonumber\\
&+\sum_{u=1}^{\ell}\sum_{l^u\in
  G_l^u}\frac{(-1)^{\ell-u+1}}{N^{\ell-u}}\sum_{r=1}^{\infty}\sum_{p=1}^{3^r}\sum_{\eta\in
C^{\epsilon^{l^u}R_{l^u}}_{r,p}}\frac{\beta^{2r}\rho^{\epsilon^{l^u}R_{l^u}}(\eta)}{r!}\nu_0(\epsilon^{l^u}R_{l^u}S_{\eta}^-)\nonumber
\end{align}
Now, we use the fact that $S_{\eta}^-=\epsilon^{\eta}R_{\eta}^-$ and we use again the
relation (\ref{E:R-to-R^-}) for each term $\epsilon^{l^u}R_{l^u}$ then
\begin{align}
&\nu_0(\epsilon^lR_l^-)=\frac{(-1)^{\ell}}{N^{\ell}}+\sum_{u=1}^{\ell}\sum_{l^u\in
G_l^u}\frac{(-1)^{\ell-u}}{N^{\ell-u}}\nu(\epsilon^{l^u}R_{l^u})+
\sum_{u=1}^{\ell}\sum_{l^u\in
  G_l^u}\frac{(-1)^{\ell-u+1}}{N^{\ell-u}}\Bigg(\nonumber\\&\sum_{r=1}^{\infty}\sum_{p=1}^{3^r}\sum_{\eta\in
C^{\epsilon^{l^u}R_{l^u}}_{r,p}}\frac{\beta^{2r}\rho^{\epsilon^{l^u}R_{l^u}}(\eta)}{r!}\frac{1}{N^u}\nu_0(\epsilon^{\eta}R_{\eta}^-)\Bigg)+\sum_{u=1}^{\ell}\sum_{l^u\in
G_l^u}\sum_{r=1}^{\infty}\sum_{p=1}^{3^r}\sum_{\eta\in C^{\epsilon^{l^u}R_{l_u}}_{r,p}}\sum_{v=1}^u\nonumber\\&\sum_{l^{u,v}\in
G_{l^u}^v}\Bigg(\frac{(-1)^{\ell-u+1}}{N^{\ell-v}}\frac{\beta^{2r}\rho^{\epsilon^{l^u}R_{l^u}}(\eta)}{r!}\nu_0(\epsilon^{l^{u,v}}R_{l^{u,v}}^-\epsilon^{\eta}R_{\eta}^-)\Big)\Bigg).\nonumber
\end{align}
the conditions $\epsilon^{\eta}=1$ in the third term and
$\epsilon^{l^{u,v}}\epsilon^{\eta}=1$ in the last term are obtained by one
application of the Proposition \ref{P:media-nula}. The result is proved by
the fact that
\begin{align}
\nu(\epsilon^lR_l)=\nu(R_l)=\sum_{i=1}^N\frac{1}{N}\nu(\sigma_i^{l_{\ell}}\sigma_i^{l_{\ell}^{\prime}}R_{\hat
l})=\nu(\epsilon^{l_{\ell}}\epsilon^{l_{\ell}^{\prime}}R_{\hat
l})=\nu(\epsilon^{\hat l}R_{\hat l})\nonumber
\end{align}
where we have used in the last two equalities the property of symmetry of
the sites and that $1=\epsilon^l=\epsilon^{\hat
  l}\epsilon^{l_{\ell}}\epsilon^{l_{\ell}^{\prime}}$. 
\qed
\end{proof}

Indeed, we can simplify more this last expression. We obtain that
\begin{proposition}\label{P:nu_0} Let
  $l=(l_1,l_1^{\prime},\ldots,l_{\ell},l_{\ell}^{\prime})$ be a sequence and
  we consider $\hat l=(l_1,\ldots,l_{\ell-1}^{\prime})$ such that
  $\epsilon^{l}=1$ and $l_i<l_i^{\prime}$ for $i=1,\ldots,\ell$. Then, we have
\begin{align}\label{E:nu_0}
&\nu_0(\epsilon^lR_l^-)=\sum_{u=1}^{\ell-2}\sum_{l^u\in
  G_l^u}\sum_{v=1}^u\sum_{l^{u,v}\in G_{l^u}^v}\frac{(-1)^{\ell-u}}{N^{\ell-v}}
\nu_0(\epsilon^{l^{u,v}}R_{l^{u,v}}^-)+\frac{1}{N^{\ell-1}}\nonumber\\&
+\sum_{u=1}^{\ell-2}\sum_{\hat l^u\in G_{\hat l}^u}\frac{\nu_0(\epsilon^{\hat
  l^u}R_{\hat l^u}^-)}{N^{\ell-1-u}}+\sum_{u=1}^{\ell-1}\sum_{\hat l^u\in
  G_{\hat
  l}^u}\sum_{r=1}^{\infty}\sum_{p=1}^{3^r}\sum_{\begin{array}{c}\eta\in
  C_{r,p}^{\epsilon^{\hat l^u}R_{\hat l^u}}\\\epsilon^{\hat l^u}\epsilon^{\eta}=1
\end{array}}\nonumber\\&\frac{\beta^{2r}\rho^{\epsilon^{\hat l^u}R_{\hat l^u}}(\eta)}{N^{\ell-1-u}r!}\nu_0(R_{\hat
  l^u}^-R_{\eta}^-)-\sum_{r=1}^{\infty}\sum_{p=1}^{3^r}\sum_{\begin{array}{c}\eta\in
  C_{r,p}^{\epsilon^lR_l}\\\epsilon^{\eta}=1\end{array}}\frac{\beta^{2r}\rho^{\epsilon^{l}R_{l}}(\eta)}{r!N^{\ell}}\nu_0(\epsilon^{\eta}R_{\eta}^-)\nonumber\\&-\sum_{r=1}^{\infty}\sum_{p=1}^{3^r}\sum_{\eta\in
  C_{r,p}^{\epsilon^lR_l}}\sum_{u=1}^{\ell}\sum_{\begin{array}{c}l^u\in G_l^u\\\epsilon^{l^u}\epsilon^{\eta}=1
\end{array}}\frac{\beta^{2r}\rho^{\epsilon^lR_l}(\eta)}{r!N^{\ell-u}}\nu_0(R_{l^u}^-R_{\eta}^-).
\end{align}
\end{proposition}
\begin{proof}
Using Proposition~\ref{P:nu_0-R-} and Equation~(\ref{E:taylor-f}) we obtain
that
\begin{align}
&\nu_0(\epsilon^lR_l^-)=\frac{(-1)^{\ell}}{N^{\ell}}+\sum_{u=1}^{\ell-1}\sum_{l^u\in
G_l^u}\frac{(-1)^{\ell-u}}{N^{\ell-u}}\nu_0(\epsilon^{l^u}R_{l^u})+\sum_{u=1}^{\ell-1}\sum_{l^u\in
G_l^u}\sum_{r=1}^{\infty}\sum_{p=1}^{3^r}\nonumber\\
&\sum_{\eta\in
  C_{r,p}^{\epsilon^{l^u}R_{l^u}}}\frac{(-1)^{\ell-u}}{N^{\ell-u}}\frac{\beta^{2r}\rho^{\epsilon^{l^u}R_{l^u}}(\eta)}{r!}\nu_0(\epsilon^{l^u}R_{l^u}S_{\eta}^-)+\nu(\epsilon^{\hat l}R_{\hat l})+\nonumber\\
&\sum_{u=1}^{\ell}\sum_{l^u\in
  G_{l}^u}\sum_{r=1}^{\infty}\sum_{p=1}^{3^r}\sum_{\begin{array}{c}\eta\in
    C_{r,p}^{\epsilon^{l^u}R_{l^u}}\\\epsilon^{\eta}=1\end{array}}\frac{(-1)^{\ell-u+1}\beta^{2r}\rho^{\epsilon^{l^u}R_{l^u}}(\eta)}{r!N^{\ell}}\nu_0(\epsilon^{\eta}R_{\eta}^-)+
\sum_{u=1}^{\ell}\nonumber\\&\sum_{l^u\in
  G_l^u}\sum_{r=1}^{\infty}\sum_{p=1}^{3^r}\sum_{\eta\in
  C_{r,p}^{\epsilon^{l^u}R_{l^u}}}\sum_{v=1}^u\sum_{\begin{array}{c}l^{u,v}\in G_{l^u}^v\\\epsilon^{l^{u,v}}\epsilon^{\eta}=1\end{array}}\frac{(-1)^{\ell-u+1}\beta^{2r}\rho^{\epsilon^{l^u}R_{l^u}}(\eta)}{r!N^{\ell-v}}\nu_0(R_{l^{u,v}}^-R_{\eta}^-).\nonumber
\end{align} 
At this point, we aplly relation~(\ref{E:R-to-R^-}) then we have
\begin{align}\label{E:inter}
&\nu_0(\epsilon^lR_l^-)=\frac{(-1)^{\ell}}{N^{\ell}}+\sum_{u=1}^{\ell-1}\sum_{l^u\in
G_l^u}\frac{(-1)^{\ell-u}}{N^{\ell-u}}\Big(\frac{1}{N^u}+\sum_{v=1}^u\sum_{l^{u,v}\in
G_{l^u}^v}\frac{1}{N^{u-v}}\times\nonumber\\&\times\nu_0(\epsilon^{l^{u,v}}R_{l^{u,v}}^-)\Big)+\sum_{u=1}^{\ell-1}\sum_{l^u\in
G_l^u}\sum_{r=1}^{\infty}\sum_{p=1}^{3^r}\sum_{\eta\in
C_{r,p}^{\epsilon^{l^u}R_{l^u}}}\frac{(-1)^{\ell-u}\beta^{2r}\rho^{\epsilon^{l^u}R_{l^u}}(\eta)}{N^{\ell-u}r!}\frac{1}{N^u}\times\nonumber\\&\times\nu_0(S_{\eta}^-)+\sum_{u=1}^{\ell-1}\sum_{l^u\in
G_l^u}\sum_{r=1}^{\infty}\sum_{p=1}^{3^r}\sum_{\eta\in C_{r,p}^{\epsilon^{l^u}R_{l^u}}}\sum_{v=1}^u\sum_{l^{u,v}\in
G_{l^u}^v}\frac{(-1)^{\ell-u}\beta^{2r}\rho^{\epsilon^{l^u}R_{l^u}}(\eta)}{N^{\ell-u}r!}\times\nonumber\\&\times\frac{1}{N^{u-v}}\nu_0(\epsilon^{l^{u,v}}R_{l^{u,v}}^-\epsilon^{\eta}R_{\eta}^-)+\nu(\epsilon^{\hat
l}R_{\hat l})+\sum_{u=1}^{\ell}\sum_{l^u\in
  G_l^u}\sum_{r=1}^{\infty}\sum_{p=1}^{3^r}\sum_{\begin{array}{c}\eta\in
    C_{r,p}^{\epsilon^{l^u}R_{l^u}}\\\epsilon^{\eta}=1\end{array}}\nonumber\\&\frac{(-1)^{\ell-u+1}\beta^{2r}\rho^{\epsilon^{l^u}R_{l^u}}(\eta)}{r!N^{\ell}}\nu_0(\epsilon^{\eta}R_{\eta}^-)+\sum_{u=1}^{\ell}\sum_{l^u\in
  G_l^u}\sum_{r=1}^{\infty}\sum_{p=1}^{3^r}\sum_{\eta\in
  C_{r,p}^{\epsilon^{l^u}R_{l^u}}}\sum_{v=1}^u\nonumber\\&\sum_{\begin{array}{c}l^{u,v}\in
  G_{l^u}^v\\\epsilon^{l^{u,v}}\epsilon^{\eta}=1\end{array}}\frac{(-1)^{\ell-u+1}\beta^{2r}\rho^{\epsilon^{l^u}R_{l^u}}(\eta)}{r!N^{\ell-v}}\nu_0(R_{l^{u,v}}^-R_{\eta}^-).
\end{align}
Now, we make some remarks. The first one
$$\frac{(-1)^{\ell}}{N^{\ell}}+\sum_{u=1}^{\ell-1}\sum_{l^u\in
  G_l^u}\frac{(-1)^{\ell-u}}{N^{\ell}}=0$$
by the binomial expansion of $(1-1)^{\ell}$. In the Equation~(\ref{E:inter}) we can
  make the assumptioms that in the second line $\epsilon^{\eta}=1$ and in the third line $\epsilon^{l^{u,v}}\epsilon^{\eta}=1$ thanks to
  Proposition~\ref{P:media-nula}. Consequently, if we do the simplifications we obtain
\begin{align}
&\nu_0(\epsilon^lR_l^-)=\sum_{u=1}^{\ell-1}\sum_{l^u\in
  G_l^u}\sum_{v=1}^u\sum_{l^{u,v}\in
  G_{l^u}^v}\frac{(-1)^{\ell-u}}{N^{\ell-v}}\nu_0(\epsilon^{l^{u,v}}R_{l^{u,v}}^-)+\nu(\epsilon^{\hat
  l}R_{\hat l})\nonumber\\&
-\sum_{r=1}^{\infty}\sum_{p=1}^{3^r}\sum_{\begin{array}{c}\eta\in
  C_{r,p}^{\epsilon^lR_{l}}\\\epsilon^{\eta}=1\end{array}}\frac{\beta^{2r}\rho^{\epsilon^lR_l}(\eta)}{r!N^{\ell}}\nu_0(\epsilon^{\eta}R_{\eta}^-)\nonumber\\
&-\sum_{r=1}^{\infty}\sum_{p=1}^{3^r}\sum_{\eta\in
  C_{r,p}^{\epsilon^lR_l}}\sum_{u=1}^{\ell}\sum_{\begin{array}{c}l^u\in G_l^u\\\epsilon^{l^u}\epsilon^{\eta}=1\end{array}}\frac{\beta^{2r}\rho^{\epsilon^lR_l}(\eta)}{r!N^{\ell-u}}\nu_0(R_{l^u}^-R_{\eta}^-).\nonumber
\end{align}
We observe that, in the first summation, the terms with $u=\ell-1$ vanish
  because $\epsilon^{l^{\ell-1}}\neq 1$ as $\epsilon^l=1$. In order to have the result annouced we look at the term $\nu(\epsilon^{\hat
  l}R_{\hat l})$ and we aplly the Proposition~\ref{P:nu-R}. Thus,
\begin{align}\label{E:inter2}
\nu(\epsilon^{\hat l}R_{\hat
  l})&=\frac{1}{N^{\ell-1}}+\sum_{u=1}^{\ell-1}\sum_{\hat l^u\in G_{\hat
  l}^u}\frac{\nu_0(\epsilon^{\hat l^u}R_{\hat
  l^u}^-)}{N^{\ell-1-u}}\nonumber\\
&+\sum_{u=1}^{\ell-1}\sum_{\hat l^u\in G_{\hat
  l}^u}\sum_{r=1}^{\infty}\sum_{p=1}^{3^r}\sum_{\begin{array}{c}\eta\in
  C_{r,p}^{\epsilon^{\hat l^u}R_{\hat l^u}}\\\epsilon^{\hat
  l^u}\epsilon^{\eta}=1\end{array}}\frac{1}{N^{\ell-1-u}}\frac{\beta^{2r}\rho^{\epsilon^{\hat
  l^u}R_{\hat l^u}}(\eta)}{r!}\nu_0(R_{\hat
  l^u}^-R_{\eta})
\end{align}
and we finish the proof observing that the term $\nu_0(\epsilon^{\hat
  l}R_{\hat l})$ is equal to zero because $\epsilon^{\hat
  l}=\epsilon_{l_{\ell}}\epsilon_{l_{\ell}^{\prime}}$ and we can apply
  Proposition~\ref{P:media-nula}. We remark that the first summation, in the Equation~(\ref{E:inter2}), vanishes if $\ell$ is equal to 2\qed
\end{proof}

Soon, we will need the following estimations 
\begin{proposition}\label{P:estimation} For all values of $\beta$ less than one, given $\ell\geq
  1$ and a sequence $l=(l_1,l_1^{\prime},\ldots,l_{\ell},l_{\ell}^{\prime})$
  we have
\begin{itemize}
\item[(a)]$\nu(S_l^-),\nu(S_l)$ and $\nu_0(S_l^-)$ are all of them $\mathrm O\left(\frac{1}{N^{\frac{\ell}{2}}}\right)$,
\item[(b)] for any function $f$ defined on $\Sigma_N^n$ and any positive
  integer $u$ and a real number $t$ in the interval $[0,1]$. There exists a
  constant $\kappa=\kappa(\beta,u,n)$ such that
$$\left|\nu_t^{(u)}(f)\right|\leq
\frac{\kappa}{N^{\frac{u}{2}}}\nu(f^2)^{\frac{1}{2}},$$
\item[(c)] for any real number $t$ in the inteval $[0,1]$ and any positive
  integer $u$: $\nu_{t}^{(u)}(S_l^-)$ and $\nu_t^{(u)}(S_l)$ are of order $\mathrm
  O\left(\frac{1}{N^{\frac{u+\ell}{2}}}\right).$
\end{itemize}
\end{proposition}
\begin{proof}
The proof of (a) can be found in \cite{BMRT} of (b) in \cite{Ta}. The result
(c)
is a consequence of (a) and (b).
\end{proof}

Given a sequence $l=(l_1,\ldots,l_{\ell}^{\prime})$ such that $\epsilon^l=1$
and $l_i<l_i^{\prime}$ for all $i=1,\ldots,\ell$ with $\ell$ a positive
integer less or equal to $2k$ where $k$ is a positive integer, we
want to find out the
expansion:
\begin{align}
\nu_0(\epsilon^lR_l^-)=\sum_{j=1}^k\frac{\lambda_j^{l}(\beta)}{N^j}+\mathrm O\left(\frac{1}{N^{k+\frac{1}{2}}}\right).\nonumber
\end{align}

First, we observe the case where $k$ is equal to one. There is no sequence $l$
such that $\ell=1$ and $\epsilon^l=1$ with $l_1<l_1^{\prime}$. Hence, we look at the following particular result

\begin{proposition}\label{P:R_{1,2}}Let $\beta$ be a positive real number less than one then
\begin{align}
\nu_0(R_{1,2}^-R_{1,2}^-)=\frac{1}{N(1-\beta^2)}+\mathrm O\left(\frac{1}{N^{\frac{3}{2}}}\right).\nonumber
\end{align}
\end{proposition}
\begin{proof}
When we use the Proposition \ref{P:nu_0}, we obtain that
\begin{align}\label{E:termo}
\nu_0(R_{1,2}^-R_{1,2}^-)&=\frac{1}{N}+\beta^2\nu_0(R_{1,2}^-R_{1,2}^-)\nonumber\\&+\sum_{r=2}^{\infty}\sum_{p=1}^{3^r}\sum_{\begin{array}{c}\eta\in
  C_{r,p}^{\epsilon_1\epsilon_2R_{1,2}}\\\epsilon^{\eta}\epsilon_1\epsilon_2=1\end{array}}\frac{\beta^{2r}\rho^{\epsilon_1\epsilon_2R_{1,2}}(\eta)}{r!}\nu_0(R_{1,2}^-R_{\eta}^-)\nonumber\\
&-\sum_{r=1}^{\infty}\sum_{p=1}^{3^r}\sum_{\begin{array}{c}\eta\in
  C_{r,p}^{R_{1,2}^2}\\\epsilon^{\eta}=1\end{array}}\frac{\beta^{2r}\rho^{R_{1,2}^2}(\eta)}{r!N^2}\nu_0(\epsilon^{\eta}R_{\eta}^-)+\nonumber\\&-\sum_{r=1}^{\infty}\sum_{p=1}^{3^r}\sum_{\eta\in
  C_{r,p}^{R_{1,2}^2}}\sum_{u=1}^2\sum_{\begin{array}{c}l^u\in
  G_{1,2,1,2}^u\\\epsilon^{l^u}\epsilon^{\eta}=1\end{array}}\frac{\beta^{2r}\rho^{R_{1,2}^2}(\eta)}{r!N^{2-u}}\nu_0(R_{l^u}^-R^-_{\eta})
\end{align}
We observe that for the first summation we can forget the condition
$\epsilon_1\epsilon_2\epsilon^{\eta}=1$ because the other elements vanish and
we can replace $C_{r,p}^{\epsilon_1\epsilon_2R_{1,2}}$ by
$C_{r,p}^{\epsilon_1\epsilon_2R_{1,2}^-}$ and
$\rho^{\epsilon_1\epsilon_2R_{1,2}}(\eta)$ by $\rho^{\epsilon_1\epsilon_2R_{1,2}^-}(\eta)$ as the functions $R_{1,2}$ and
$R_{1,2}^-$ are different by only one constant that does not depends on any
configuration which implies:
\begin{align}\label{E:termo1}
\sum_{r=2}^{\infty}\sum_{p=1}^{3^r}&\sum_{\begin{array}{c}\eta\in
    C^{\epsilon_1\epsilon_2R_{1,2}}_{r,p}\\\epsilon_1\epsilon_2\epsilon^{\eta}=1\end{array}}\frac{\beta^{2r}\rho^{\epsilon_1\epsilon_2R_{1,2}}(\eta)}{r!}\nu_0(R_{1,2}^-R^-_{\eta})=\nonumber\\&=\sum_{r=2}^{\infty}\sum_{p=1}^{3^r}\sum_{\eta\in
C_{r,p}^{\epsilon_1\epsilon_2R_{1,2}^-}}\frac{\beta^{2r}\rho^{\epsilon_1\epsilon_2R_{1,2}^-}(\eta)}{r!}\nu_0(\epsilon_1\epsilon_2R^-_{1,2}\epsilon^{\eta}R_{\eta}^-)\nonumber\\&=\nu_{t_1}^{(2)}(\epsilon_1\epsilon_2R_{1,2}^-)=\mathrm O\left(\frac{1}{N^{\frac{3}{2}}}\right)
\end{align}
In the last two equality, we have just used the fact that the summation is the
error term in the Taylor expansion of order one for the function
$\epsilon_1\epsilon_2R^-_{1,2}$ and the Proposition~\ref{P:estimation} item
(c) where $t_1$ is a real number in the interval $[0,1]$.

Now, we dont consider the constraints $\epsilon^{\eta}=1$ and
$\epsilon^{l^u}\epsilon^{\eta}=1$ in the other summation then using relation~(\ref{E:R-to-R^-}) the other summations can be rewritten as
\begin{align}\label{E:termo2}
&\sum_{r=1}^{\infty}\sum^{3^r}_{p=1}\sum_{\eta\in
    C^{R_{1,2}^2}_{r,p}}\frac{\beta^{2r}\rho^{R_{1,2}^2}(\eta)}{r!}\Bigg(\frac{\nu_0(\epsilon^{\eta}R_{\eta}^-)}{N^2}+\sum_{u=1}^2\sum_{l^u\in
    G_{1,2,1,2}^u}\frac{1}{N^{2-u}}\times\nonumber\\&\times\nu_0(\epsilon^{l^u}R_{l^u}^-\epsilon^{\eta}R_{\eta}^-)\Bigg)=\sum_{r=1}^{\infty}\sum_{p=1}^{3^r}\sum_{\eta\in
    C_{r,p}^{R_{1,2}^2}}\frac{\beta^{2r}\rho^{R_{1,2}^2}(\eta)}{r!}\nu_0(R_{1,2}^2\epsilon^{\eta}R_{\eta}^-)\nonumber\\&=\nu_{t_2}^{(1)}(R_{1,2}^2)=\mathrm
    O\left(\frac{1}{N^{\frac{3}{2}}}\right).
\end{align}
where the last summation can be seen as the error term in the Taylor expansion of
order zero for
the function $\nu(R_{1,2}^2)$ and we used again the Proposition~\ref{P:estimation}
item (c).

Putting together the Equations~(\ref{E:termo}),(\ref{E:termo1}) and
(\ref{E:termo2}) we obtain
\begin{align}
\nu_0(R_{1,2}^-R_{1,2}^-)&=\frac{1}{N}+\beta^2\nu_0(R_{1,2}^-R_{1,2}^-)+\mathrm O\left(\frac{1}{N^{\frac{3}{2}}}\right)\nonumber
\end{align}
which concludes the proof \qed
\end{proof}
\section{The proof of Lemma~\ref{L:main}}
\begin{proof}
We observe that
$\nu_0(\epsilon^lR_l^-)=\mathrm O\left(\frac{1}{N^{\frac{\ell}{2}}}\right)$ by
Proposition~\ref{P:estimation} item (a) which implies $\lambda_j^l=0$ for any
nonpositive integer number $j$ and any sequence $l$.

\textit{(Step 0): The case $k=1$}

If $k=1$ the sequences $l=(l_1,l_1^{\prime},\ldots,l_{\ell},l_{\ell}^{\prime})$ such that $\ell\leq 2$ and
$\epsilon^{\ell}=1$ are of the type $(l_1,l_1^{\prime},l_1,l_1^{\prime})$. As
$R_{l_1,l_1^{\prime}}^-$ depends on only two configurations we have that
\begin{align}
\nu_0(R_{l_1,l_1^{\prime}}^-R_{l_1,l_1^{\prime}}^-)=\nu_0(R_{1,2}^-R_{1,2}^-)=\frac{1}{N(1-\beta^2)}+\mathrm O\left(\frac{1}{N^{\frac{3}{2}}}\right),
\end{align}
where we have used in the last equality the
Proposition~\ref{P:R_{1,2}}. Hence, we have proved (a).

\textit{(Step 1): the case $k\geq 2$ with $\ell=2k$}

The proof is by induction on $k$. The case $k=1$ has been proved in the
\textit{Step 0}. Now, we suppose to know all of the coefficients $\lambda_j^{\ell}$
for all sequences $l$ such that $\ell\leq 2k-2$ and $j=1,\ldots,k-1$. First,
we evaluate the coefficient $\lambda_k^l$ for a sequence $l$ such that
$\ell=2k$. We apply the hypothesis of induction in the respectively terms in
the expression of Proposition~\ref{P:nu_0} which implies
\begin{align}
\nu_0(\epsilon^lR_l^-)=A_1+A_2+A_3+A_4\nonumber 
\end{align}
The term of the first and second summations in (\ref{E:nu_0}):
\begin{align}
A_1&=\sum_{u=1}^{2k-2}\sum_{l^u\in
  G_l^u}\sum_{v=1}^u\sum_{l^{u,v}\in
  G_{l^u}^v}\frac{(-1)^{2k-u}}{N^{2k-v}}\left(\sum_{j=1}^{k-1}\frac{\lambda_j^{l^{u,v}}}{N^j}+\mathrm
  O\left(\frac{1}{N^{k-\frac{1}{2}}}\right)\right)+\nonumber\\&+\frac{1}{N^{2k-1}}+\sum_{u=1}^{2k-2}\sum_{\hat l^u\in G_{\hat
  l}^u}\frac{1}{N^{2k-1-u}}\left(\sum_{j=1}^{k-1}\frac{\lambda_j^{\hat
  l^u}}{N^j}+\mathrm
  O\left(\frac{1}{N^{k-\frac{1}{2}}}\right)\right),\nonumber
\end{align}
the term of the third summation in (\ref{E:nu_0}) with $u\leq 2k-2$:
\begin{align}
A_2&=\sum_{u=1}^{2k-2}\sum_{\hat
  l^u\in G_{\hat
  l}^u}\Bigg[\sum_{r=1}^{2k-2-u}\sum_{p=1}^{3^r}\sum_{\begin{array}{c}\eta\in
  C_{r,p}^{\epsilon^{\hat l^u}R_{\hat l^u}}\\\epsilon^{\hat
  l^u}\epsilon^{\eta}=1\end{array}}\frac{\beta^{2r}\rho^{\epsilon^{\hat
  l^u}R_{\hat l^u}}(\eta)}{N^{2k-1-u}r!}\Bigg(\sum_{j=1}^{k-1}\frac{\lambda_j^{(\hat
  l^u,\eta)}}{N^j}\nonumber\\&+\mathrm
  O\left(\frac{1}{N^{k-\frac{1}{2}}}\right)\Bigg)+\sum_{r\geq 2k-1-u}\sum_{p=1}^{3^r}\sum_{\begin{array}{c}\eta\in
  C_{r,p}^{\epsilon^{\hat l^u}R_{\hat l^u}}\\\epsilon^{\hat
  l^u}\epsilon^{\eta}=1\end{array}}\frac{\beta^{2r}\rho^{\epsilon^{\hat
  l^u}R_{\hat l^u}}(\eta)\nu_0(R_{\hat
  l^u}^-R_{\eta}^-)}{r!N^{2k-1-u}}\Bigg].\nonumber
\end{align}
Now, the term with $u=2k-1$, we remark that for $r=1$ it appears the term
$\epsilon^lR_l$. Thus,
\begin{align}
A_3&=\beta^{2}\nu_0(\epsilon^{l}R_{l}^-)+\sum_{r\geq
  2}\sum_{p=1}^{3^r}\sum_{\begin{array}{c}\eta\in C_{r,p}^{\epsilon^{\hat
  l}R_{\hat l}}\\\epsilon^{\hat
  l}\epsilon^{\eta}=1\end{array}}\frac{\beta^{2r}\rho^{\epsilon^{\hat
  l}R_{\hat l}}(\eta)}{r!}\nu_0(R_{\hat
  l}^-R_{\eta}^-),\nonumber
\end{align}
the last terms:
\begin{align}
A_4=-\sum_{r=1}^{\infty}\sum_{p=1}^{3^r}\sum_{\begin{array}{c}\eta\in
  C_{r,p}^{\epsilon^lR_l}\end{array}}\frac{\beta^{2r}\rho^{\epsilon^lR_l}(\eta)}{r!}
&\Bigg[\frac{1}{N^{2k}}\nu_0(\epsilon^{\eta}R_{\eta}^-)\nonumber\\&+\sum_{u=1}^{2k}\sum_{l^u\in
  G_l^u}\frac{1}{N^{2k-u}}\nu_0(\epsilon^{l^u}R_{l^u}^-\epsilon^{\eta}R_{\eta}^-)\Bigg].\nonumber
\end{align}
In the term $A_4$ we removed the constraints $\epsilon^{\eta}=1$ and
$\epsilon^{l^u}\epsilon^{\eta}=1$ as for the other $\eta$ these
terms vanish. The result follows if we look at the coefficient
of order $\frac{1}{N^k}$ for each $A_j$ with $j=1,\ldots,4$ and if we demonstrate that the other terms are of
order $\mathrm O\left(\frac{1}{N^{k+\frac{1}{2}}}\right)$. Indeed, each term
\begin{align}
\sum_{r\geq 2k-1-u}\sum_{p=1}^{3^r}\sum_{\begin{array}{c}\eta\in
    C_{r,p}^{\epsilon^{\hat l^u}R_{\hat l_ u}}\\\epsilon^{\hat
    l^u}\epsilon^{\eta}=1\end{array}}&\frac{\beta^{2r}\rho^{\epsilon^{\hat
    l^u}R_{\hat l^u}}(\eta)\nu_0(R_{\hat
    l^u}^-R_{\eta}^-)}{r!N^{2k-1-u}}=\nonumber\\&=\frac{1}{N^{2k-1-u}}\nu_t^{(2k-1-u)}(\epsilon^{\hat l^u}R_{\hat l^u}),\nonumber
\end{align}
that appears in $A_2$, satisfies for a positive real number $t$ in the interval $[0,1]$ and we could
remove the condition $\epsilon^{\hat l^u}\epsilon^{\eta}=1$ because for other
values of $\eta$ it vanishes and we could replace $\eta\in \epsilon^{\hat
  l^u}R_{\hat l^u}$ by $\eta\in\epsilon^{\hat l^u}R_{\hat l^u}^-$ because these
functions depends on the same configurations and we have used the fact that it
could be seen as a error term in the expansion of $\nu(\epsilon^{\hat
  l^u}R_{\hat l^u}^-)$ of order $2k-2-u$. As $u\leq 2k-2$ and using
Propostion~{\ref{P:estimation}} item (c) each one of these terms are $\mathrm
O\left(\frac{1}{N^{k+\frac{1}{2}}}\right)$. By the same thought for a real
number $t_1$ in the interval $[0,1]$, we can show
that, the terms which appear in $A_3$,
\begin{align}
\sum_{r\geq 2}\sum_{p=1}^{3^r}\sum_{\begin{array}{c}\eta\in
    C_{r,p}^{\epsilon^{\hat l}R_{\hat l}}\\\epsilon^{\hat
    l}\epsilon^{\eta}=1\end{array}}\frac{\beta^{2r}\rho^{\epsilon^{\hat
    l}R_{\hat l}}(\eta)}{r!}\nu_0(R_{\hat
    l}^-R_{\eta}^-)=\nu_{t_1}^{(2)}(\epsilon^{\hat l}R_{\hat l}^-)=\mathrm O\left(\frac{1}{N^{k+\frac{1}{2}}}\right).
\end{align}
To finish, we look at the terms in $A_4$ and we use the relation~(\ref{E:R-to-R^-}) then
\begin{align}
&\sum_{r=1}^{\infty}\sum_{p=1}^{3^r}\sum_{\eta\in
  C_{r,p}^{\epsilon^lR_l}}\frac{\beta^{2r}\rho^{\epsilon^lR_l}(\eta)}{r!}\Bigg[\frac{1}{N^{2k}}\nu_0(\epsilon^{\eta}R_{\eta}^-)+\sum_{u=1}^{2k}\sum_{l^u\in
G_{l}^u}\frac{1}{N^{2k-u}}\times\nonumber\\&\times\nu_0(\epsilon^{l^u}R_{l^u}^-\epsilon^{\eta}R_{\eta}^-)\Bigg]=\sum_{r\geq
  1}\sum_{p=1}^{3^r}\sum_{\eta\in
  C_{r,p}^{\epsilon^{l}R_l}}\frac{\beta^{2r}\rho^{\epsilon^lR_l}(\eta)}{r!}\nu_0(\epsilon^lR_l\epsilon^{\eta}R_{\eta}^-)\nonumber\\&=\nu_{t_2}^{(1)}(\epsilon^lR_l)=\mathrm O\left(\frac{1}{N^{k+\frac{1}{2}}}\right)
\end{align}
where we have used the fact that the last summation can be seen as a error
term in the Taylor's expansion of order zero for the function
$\nu(\epsilon^lR_l)$. Afterwards, we applied again the
Proposition~{\ref{P:estimation}} item (c) which proves (b).

\textit{(Step 2:) the case $k\geq 2$ with $\ell<2k$}

Now, we go on with the induction and we suppose the assumptions of item
(c). Then, we apply Proposition~\ref{P:nu_0} for a sequence $l$ such that
$2\leq \ell=\hat k\leq 2k$ which produces
\begin{align}
\nu_0(\epsilon^lR_l^-)=B_1+B_2+B_3+B_4
\end{align}
where $B_1$ is associated to the two first summations in (\ref{E:nu_0}) and we have applied the
hypothesis of induction claimed in the item (b):
\begin{align}
B_1&=\sum_{u=1}^{\hat k-1}\sum_{l^u\in
  G_l^u}\sum_{v=1}^u\sum_{l^{u,v}\in G_{l^u}^v}\frac{(-1)^{\hat k-u}}{N^{\hat
  k-v}}\left(\sum_{j=1}^{k-1}\frac{\lambda_j^{l^{u,v}}}{N^j}+\mathrm
  O\left(\frac{1}{N^{k-\frac{1}{2}}}\right)\right)+\nonumber\\&+\frac{1}{N^{\hat
  k-1}}+\sum_{u=1}^{\hat k-2}\sum_{\hat l^u\in G_{\hat l}^u}\frac{1}{N^{\hat
  k-1-u}}\left(\sum_{j=1}^{k-1}\frac{\lambda_j^{\hat l^u}}{N^j}+\mathrm
  O\left(\frac{1}{N^{k-\frac{1}{2}}}\right)\right).\nonumber
\end{align}
$B_2$ has $u\leq \hat k-2$ in the third summation at (\ref{E:nu_0}) and we apply again the
hypothesis of item (b):
\begin{align}
B_2&=\sum_{u=1}^{\hat
  k-2}\sum_{\hat l^u\in G_{\hat
  l}^u}\Bigg[\sum_{r=1}^{2k-2-u}\sum_{p=1}^{3^r}\sum_{\begin{array}{c}\eta\in
  C_{r,p}^{\epsilon^{\hat l}R_{\hat l}}\\\epsilon^{\hat
  l^u}\epsilon^{\eta}=1\end{array}}\frac{\beta^{2r}\rho^{\epsilon^{\hat
  l}R_{\hat l}}(\eta)}{r!N^{\hat
  k-1-u}}\Bigg(\sum_{j=1}^{k-1}\frac{\lambda_j^{(\hat l^u,\eta)}}{N^j}+\nonumber\\&+\mathrm
  O\left(\frac{1}{N^{k-\frac{1}{2}}}\right)\Bigg)+\sum_{r\geq
  2k-1-u}\sum_{p=1}^{3^r}\sum_{\begin{array}{c}\eta\in C_{r,p}^{\epsilon^{\hat
  l^u}R_{\hat l^u}}\\\epsilon^{\hat l^u}\epsilon^{\eta}=1
  \end{array}}\frac{\beta^{2r}\rho(\eta)}{N^{\hat
  k-1-u}r!}\nu_0(\epsilon^{\hat l^u}R_{\hat
  l^u}^-\epsilon^{\eta}R_{\eta}^-)\Bigg].\nonumber
\end{align}
$B_3$ is associated to the third summation when $u=\hat k-1$. Now, we also
apply the hypothesis claimed in the item (c):
\begin{align}
B_3&=\beta^2\nu_0(\epsilon^lR_l^-)+\Bigg[\sum_{r=2}^{2k-\hat
  k+1}\sum_{p=1}^{3^r}\nonumber\\&\sum_{\begin{array}{c}\eta\in
  C_{r,p}^{\epsilon^{\hat l}R_{\hat l}}\\\epsilon^{\hat
  l}\epsilon^{\eta}=1\end{array}}\frac{\beta^{2r}\rho^{\epsilon^{\hat
  l}R_{\hat l}}(\eta)}{r!}\Bigg(\sum_{j=1}^k\frac{\lambda_j^{(\hat
  l,\eta)}}{N^j}+\mathrm
  O\left(\frac{1}{N^{k+\frac{1}{2}}}\right)\Bigg)\nonumber\\&+\sum_{r\geq
  2k-\hat k+2}\sum_{p=1}^{3^r}\sum_{\begin{array}{c}\eta\in
  C_{r,p}^{\epsilon^{\hat l}R_{\hat l}}\\\epsilon^{\hat
  l}\epsilon^{\eta}=1\end{array}}\frac{\beta^{2r}\rho^{\epsilon^{\hat
  l}R_{\hat l}}(\eta)}{r!}\nu_0(\epsilon^{\hat
  l}R_{\hat l}^-\epsilon^{\eta}R_{\eta}^-)\Bigg],\nonumber
\end{align}
the last terms appear in $B_4$:
\begin{align}
B_4&=-\sum_{r=1}^{2k-\hat k}\sum_{p=1}^{3^r}\sum_{\begin{array}{c}\eta\in
    C_{r,p}^{\epsilon^lR_l}\\\epsilon^{\eta}=1\end{array}}\frac{\beta^{2r}\rho^{\epsilon^lR_l}(\eta)}{r!N^{\hat
    k}}\left(\sum_{j=1}^{k-1}\frac{\lambda_j^{\eta}}{N^{j}}+\mathrm
    O\left(\frac{1}{N^{k-\frac{1}{2}}}\right)\right)-\nonumber\\&
-\sum_{r\geq 2k-\hat k+1}\sum_{p=1}^{3^r}\sum_{\begin{array}{c}\eta\in
    C_{r,p}^{\epsilon^lR_l}\\\epsilon^{\eta}=1\end{array}}\frac{\beta^{2r}\rho^{\epsilon^lR_l}(\eta)}{r!N^{\hat
    k}}\nu_0(\epsilon^{\eta}R_{\eta}^-)\nonumber\\&-\sum_{r=1}^{2k-\hat
    k}\sum_{p=1}^{3^r}\sum_{\eta\in
    C_{r,p}^{\epsilon^lR_l}}\frac{\beta^{2r}\rho^{\epsilon^lR_l}(\eta)}{r!}\Bigg[\sum_{u=1}^{\hat
    k-1}\sum_{\begin{array}{c}l^u\in
    G_l^u\\\epsilon^{l^u}\epsilon^{\eta}=1\end{array}}\frac{1}{N^{\hat
    k-u}}\Bigg(\sum_{j=1}^{k-1}\frac{\lambda_j^{(l^u,\eta)}}{N^j}+\nonumber\\&+\mathrm
    O\left(\frac{1}{N^{k-\frac{1}{2}}}\right)\Bigg)+\sum_{j=1}^k\frac{\lambda_j^{(l,\eta)}}{N^j}+\mathrm
    O\left(\frac{1}{N^{k+\frac{1}{2}}}\right)\Bigg]-\nonumber\\&-\sum_{r\geq 2k-\hat
    k+1}\sum_{p=1}^{3^r}\sum_{\eta\in C_{r,p}^{\epsilon^lR_l}}\sum_{u=1}^{\hat
    k}\sum_{\begin{array}{c}l^u\in G_l^u\\\epsilon^{l^u}\epsilon^{\eta}=1\end{array}}\frac{\beta^{2r}\rho^{\epsilon^lR_l}(\eta)}{r!N^{\hat k-u}}\nu_0(\epsilon^{l^u}R_{l^u}^-\epsilon^{\eta}R_{\eta}^-).\nonumber
\end{align}
Eventually, we obtain the result by observing the coefficient of the term
$\frac{1}{N^k}$ in the expressions $B_1,B_2,B_3$ and $B_4$ and by the fact
that the other terms are $\mathrm O\left(\frac{1}{N^{k+\frac{1}{2}}}\right)$.

For a fixed sequence $\hat l^u$ each term, that appears in $B_2$,
\begin{align}
&\sum_{r\geq 2k-1-u}\sum_{p=1}^{3^r}\sum_{\begin{array}{c}\eta\in
    C_{r,p}^{\epsilon^{\hat l^u}R_{\hat l^u}}\\\epsilon^{\hat
    l^u}\epsilon^{\eta}=1\end{array}}\frac{\beta^{2r}\rho^{\epsilon^{\hat
    l^u}R_{\hat l^u}}(\eta)}{N^{\hat
    k-1-u}r!}\nu_0(\epsilon^{\hat l^u}R_{\hat
    l^u}^-\epsilon^{\eta}R_{\eta}^-)=\nonumber\\&\frac{1}{N^{\hat
    k-1-u}}\nu_{t_1}^{(2k-1-u)}(\epsilon^{\hat l^u}R_{\epsilon^{\hat
    l^u}}^-)=\frac{1}{N^{\hat k-1-u}}\mathrm
    O\left(\frac{1}{N^{k-\frac{1}{2}}}\right)=\mathrm O\left(\frac{1}{N^{k+\frac{1}{2}}}\right)
\end{align}
where we have just used the fact that this term can be seen as a error term in
the Taylor expansion of order $2k-u$ of the function $\nu(\epsilon^{\hat
  l^u}R_{\hat l^u})$ and afterwards the Proposition~\ref{P:estimation} item
(c) and the fact that $u\leq \hat k-2$. By similar arguments the term, that
appears in $B_3$,
\begin{align}
\sum_{r\geq 2k-\hat k+1}\sum_{p=1}^{3^r}\sum_{\begin{array}{c}\eta\in
    C_{r,p}^{\epsilon^{\hat l}R_{\hat l}}\\\epsilon^{\hat l}\epsilon^{\eta}=1\end{array}}&\frac{\beta^{2r}\rho^{\epsilon^{\hat
    l}R_{\hat l}}(\eta)}{r!}\nu_0(\epsilon^{\hat
    l}R_{\hat l}^-\epsilon^{\eta}R_{\eta}^-)\nonumber\\&=\nu_{t_2}^{(2k-\hat k+2)}(\epsilon^{\hat
    l}R_{\hat l}^-)=\mathrm O\left(\frac{1}{N^{k+\frac{1}{2}}}\right).\nonumber
\end{align}

To finish, we look at the terms that appear in $B_4$ which we can remove the
constraints $\epsilon^{\eta}=1$ and $\epsilon^{l^u}\epsilon^{\eta}=1$ because
the other vanish and
usign relation~(\ref{E:R-to-R^-}) we have that
\begin{align}
&\sum_{r\geq 2k-\hat k+1}\sum_{p=1}^{3^r}\sum_{\eta\in C_{r,p}^{\epsilon^l
    R_l}}\frac{\beta^{2r}\rho^{\epsilon^lR_l}(\eta)}{r!}\Bigg[\nu_0(\epsilon^{\eta}R_{\eta}^-)+\sum_{u=1}^{\hat
k}\sum_{l^u\in
    G_l^u}\frac{\nu_0(\epsilon^{l^u}R_{l^u}^-\epsilon^{\eta}R_{\eta}^-)}{N^{\hat
    k-u}} \Bigg]=\nonumber\\&=\sum_{r\geq 2k-\hat
    k+1}\sum_{p=1}^{3^r}\sum_{\eta\in
    C_{r,p}^{\epsilon^lR_l}}\frac{\beta^{2r}\rho^{\epsilon^lR_l}(\eta)}{r!}\nu_0(\epsilon^lR_l^-\epsilon^{\eta}R_{\eta}^-)=\nonumber\\&=\nu_{t_3}^{(2k-\hat
    k+1)}(\epsilon^lR_l^-)=\mathrm O\left(\frac{1}{N^{k+\frac{1}{2}}}\right).
\end{align}
where we have just used that $\epsilon^lR_l^-$ and $\epsilon^lR_l$ depends on
the same configurations. We can see the last summation as the error term of
the Taylor's expansion of order $2k-\hat k$ for $\nu(\epsilon^l R_l^-)$ and we
can make use of the Proposition~{\ref{P:estimation}} item (c) for concluding the proof\qed
\end{proof}

\section{The proof of Theorem~\ref{T:main}}
\begin{proof}
Using the property of symmetry of the sites, we obtain that
\begin{align}
\nu(R_{1,2}^2)=\sum_{i=1}^N\frac{\nu(\sigma_i^1\sigma_i^2R_{1,2})}{N}=\nu(\epsilon_1\epsilon_2R_{1,2})=\frac{1}{N}+\nu(\epsilon_1\epsilon_2R_{1,2}^-)\nonumber
\end{align}
where we have just used the relation
$R_{1,2}=\frac{\epsilon_1\epsilon_2}{N}+R_{1,2}^-$. Applying the
Equation~(\ref{E:taylor-f}), we have
\begin{align}
\nu(R_{1,2}^2)=\frac{1}{N}+\sum_{\ell=1}^{\infty}\sum_{m=1}^{3^{\ell}}\sum_{l\in
C_{\ell,m}^{\epsilon_1\epsilon_2R_{1,2}^-}}\frac{\beta^{2\ell}\rho^{\epsilon_1\epsilon_2R_{1,2}^-}(l)}{\ell!}\nu_0(\epsilon_1\epsilon_2R_{1,2}^-\epsilon^lR_l^-)
\end{align}
As the term
\begin{align}
\sum_{\ell\geq 2k_0}\sum_{m=1}^{3^{\ell}}\sum_{l\in
  C_{\ell,m}^{\epsilon_1\epsilon_2R_{1,2}^-}}&\frac{\beta^{2\ell}\rho^{\epsilon_1\epsilon_2R_{1,2}^-}(\eta)}{\ell!}\nu_0(\epsilon_1\epsilon_2R_{1,2}^-\epsilon^{l}R_l^-)=\nonumber\\&=\nu_{t_1}^{(2k_0)}(\epsilon_1\epsilon_2R_{1,2}^-)=\mathrm O\left(\frac{1}{N^{k_0+\frac{1}{2}}}\right)\nonumber
\end{align}
because it can be seen as the error term of order $2k_0$ of the function
$\nu(\epsilon_1\epsilon_2R_{1,2}^-)$ and we have used the
Proposition~\ref{P:estimation} item (c). Now, we put the constraint
$\epsilon_1\epsilon_2\epsilon^l=1$ because the other terms vanish by
Proposition~\ref{P:media-nula} and we apply the Lemma~\ref{L:main}. Thus,
\begin{align}
\nu(R_{1,2}^2)=\frac{1}{N}&+\sum_{\ell=1}^{2k_0-1}\sum_{m=1}^{3^{\ell}}\sum_{\begin{array}{c}l\in
  C_{\ell,m}^{\epsilon_1\epsilon_2R_{1,2}^-}\\\epsilon_1\epsilon_2\epsilon^l=1\end{array}}\frac{\beta^{2\ell}\rho^{\epsilon_1\epsilon_2R_{1,2}^-}(l)}{\ell!}\sum_{j=1}^{k_0}\lambda_j^{(1,2,l)}+\nonumber\\&+\mathrm O\left(\frac{1}{N^{k_0+\frac{1}{2}}}\right)
\end{align} 
which concludes the proof \qed
\end{proof}

\end{document}